\documentclass{amsart}
\usepackage{amssymb,amsmath,amsthm}

\newtheorem{thm}{Theorem}[section]
\newtheorem{lem}[thm]{Lemma}

\newtheorem{prop}[thm]{Proposition}

\theoremstyle{remark}
\newtheorem{exm}{Example}

\theoremstyle{definition}
\newtheorem{defi}[thm]{Definition}

\newcommand {\Nz} {\mathbb{N}}
\newcommand {\Zz} {\mathbb{Z}}

\newcommand {\Rz} {\mathbb{R}}
\newcommand {\Cz} {\mathbb{C}}

\newcommand {\T} {\mathbb{S}}
\newcommand {\TT} {\mathbb{T}}

\newcommand{\BAR}{\overline}

\newcommand{\ii}{\mathrm{i}}

\DeclareMathOperator{\Irr}{Irr}

\DeclareMathOperator{\Sym}{Sym}

\DeclareMathOperator{\SL}{SL}

\DeclareMathOperator{\Tau}{{\mathcal T}}

\newcommand{\abin}{a}

\newcommand{\defst}[1]{{\it #1}}

\begin{document}

\title{Fusion algebras for imprimitive complex
reflection groups}
\author{Michael Cuntz}
\address{Michael Cuntz, Universit\"at Kaiserslautern,
Postfach 3049, 67653 Kaiserslautern}
\email{cuntz@mathematik.uni-kl.de}

\begin{abstract}
We prove that
the Fourier matrices for the imprimitive complex reflection groups
introduced by Malle in \cite{MUG} define
fusion algebras with not necessarily positive
but integer structure constants. Hence they define $\Zz$-algebras.
As a result, we obtain that all known Fourier matrices belonging
to spetses define algebras with integer structure constants.
\end{abstract}

\maketitle

\section{Introduction}

In his classification of irreducible characters of a finite group of Lie type,
Lusztig develops a theory in which a so-called non abelian
Fourier transform emerges.
This is a matrix which only depends on the Weyl group of the group of Lie type.
Geck and Malle \cite{mGgM} set up a system of axioms based on the properties such a
Fourier matrix has. Using this system Brou\'e, Malle and Michel construct
analogous transformations for the spetses, which until now remain mysterious objects.

In \cite{MUG} Malle defines unipotent degrees for
the imprimitive complex reflection groups.
The transformation matrix from the fake degrees to
these unipotent degrees
defines an algebra via the formula of Verlinde.
In the present article, we show that these matrices yield algebras with
integer structure constants. We prove this by looking at exterior powers
of group rings of cyclic groups.

We start by giving a definition of the specific
type of $\Zz$-algebra we will study. It has roughly
the properties of a table algebra, though it is
not a $\Cz$-algebra and its structure constants may be
negative. The based rings Lusztig introduces in \cite{Lus1}
are also a variation of the algebras we look at.
A much more extensive investigation of such algebras
can be found in the authors dissertation \cite{mC}.

We then recall the definition of the matrices from
\cite{MUG} and explain their relation to
exterior powers.
In the following section we examine the algebras
belonging to exterior powers so that we can
prove the result in the next section.
Finally, we show a connection to the Kac-Peterson matrices
attached to affine Kac-Moody algebras:
exterior powers of matrices of type $A^{(1)}_1$ are matrices of
type $C^{(1)}_l$.

This article is a part of chapter $5$ and $6$ of \cite{mC}.
I am very grateful to my supervisor G. Malle for many
helpful discussions.

\section{Basic definitions}

\begin{defi}
Let $R$ be a finitely generated commutative $\Zz$-algebra
which is a free $\Zz$-module with basis
$B=\{b_0=1,\ldots,b_{n-1}\}$ and \defst{structure constants}
\[ b_i b_j = \sum_k N_{ij}^k b_k, \quad N_{ij}^k\in \Zz \]
for $0\le i,j < n$.
Assume that there is an involution $\sim : R \rightarrow R$
which is a $\Zz$-module homomorphism such that
\[ \tilde B= B,\quad N_{\tilde i\tilde j}^k=N_{i j}^{\tilde k},
\quad N_{\tilde ij}^0 = \delta_{i,j} \]
for all $0\le i,j,k < n$,
where $\tilde i$ is the index with $\tilde{b_i}=b_{\tilde i}$.
Then we call $(R,B)$ a \defst{$\Zz$-based ring}.
\end{defi}
Remark that if the involution $\sim$ exists, then it is
unique by the third equation above.
The second equation expresses that $\sim$ is an algebra
homomorphism.
Remark also, that if we replace an element $b\in B$ by $-b$,
then the new basis spans the same algebra,
but $\sim$ does not necessarily exist anymore
(with respect to the new basis).
\begin{exm}
Let $G$ be a finite group. Then the character ring
of $G$ is a $\Zz$-based ring with basis $\Irr(G)$
(the irreducible characters)
and non negative structure constants,
where multiplication is just tensor product.
The involution $\sim$ is complex conjugation on
the characters.
\end{exm}
\begin{exm}
The representation ring of the quantum double of a finite group
(see \cite{CGR}) is a $\Zz$-based ring where the basis is again
the set of irreducible representations.
\end{exm}
The $\Zz$-based rings are a generalization of algebras
with the properties of representation rings.
There are many other such generalizations.
One of them are the table algebras, to which the $\Zz$-based
rings with non negative structure constants belong (viewed as $\Cz$-algebras).
All $\Zz$-based rings are generalized table algebras (GT-algebras, \cite{zAeFmM}).
But GT-algebras do not have the properties which
we will need. Another structure is the so-called
$C$-algebra, which has an elaborate
structure theory \cite{hB}.
Unfortunately, the $\Zz$-based rings of the present article
are not always $C$-algebras.

If $(R,B)$ is a $\Zz$-based ring, then we have a linear map
$\tau : R\rightarrow \Cz$ defined by $\tau(b_i)=\delta_{0,i}$.
The map
\[\langle\:,\:\rangle : R\times R\rightarrow \Zz,\quad
\langle r,r' \rangle:=\tau(\tilde r r')\]
for $r,r'\in R$ behaves like an inner product with
orthonormal basis $B$ because $r = \sum_{b\in B} \langle b,r \rangle b$
for all $r\in R$. The set $\tilde B$ is the basis
dual to $B$ with respect to this inner product.
Extending $\langle\:,\:\rangle$ to the
$\Cz$-algebra $R_\Cz:=R\otimes_\Zz \Cz$, one can prove
that $R_\Cz$ is semisimple (compare \cite{mC}, 1.2).
\begin{prop}
Let $R$ be a $\Zz$-based ring.
Then the algebra $R_\Cz:=R\otimes_\Zz \Cz$ is semisimple.
\end{prop}
\begin{proof}
Extend $\sim$ and $\tau$ to $R_\Cz$:
\[ \widetilde{r \otimes z} := \tilde r \otimes \bar z,
\quad \tau' : r \otimes z \mapsto z \tau(r), \]
where $r\in R$, $z \in \Cz$.
This yields a hermitian positive definite sesquilinear form
$\langle r,r' \rangle := \tau'(\tilde r r')$.
If $\mathfrak I$ is a left ideal in $R_\Cz$, then the orthogonal complement
\[ \mathfrak I^\perp := \{ r \in R_\Cz \mid \langle r,r' \rangle = 0
 \quad \forall r' \in \mathfrak I \} \]
is a left ideal too:
\[ \langle tr,r' \rangle = \tau'(\widetilde{tr} r') =
\tau'(\tilde r \tilde t r') = \langle r,\tilde t r' \rangle = 0 \]
for all $r \in \mathfrak I^\perp$, $t \in R_\Cz$ and $r' \in \mathfrak I$.
The claim follows.
\end{proof}

Now $R_\Cz$ is a commutative semisimple algebra over an algebraically
closed field, so by the theorem of Wedderburn-Artin it is isomorphic
as a $\Cz$-algebra to $\Cz^n$ with componentwise multiplication.
By choosing $B$ as a basis for $R_\Cz$ and the canonical basis
$\{e_i\}_i$ with $e_i e_j=\delta_{i,j} e_i$ for all $i,j$ for $\Cz^n$, an
isomorphism $\varphi$ is described by a matrix $s$ which we will call
an \defst{$s$-matrix of} $(R,B)$:
\[ \varphi(b_i) = \sum_k s_{ki} e_k. \]
Remark that this matrix depends on the choice of the isomorphism
$\varphi$. Another isomorphism would differ from $\varphi$ by
a $\Cz$-algebra automorphism of $\Cz^n$,
so an $s$-matrix is unique up to a permutation of rows.

The rows of $s$ are the one-dimensional representations of
$R$ because $s_{ki}s_{kj}=\sum_l N_{ij}^l s_{kl}$ for all $k,i,j$.
They are orthogonal (see \cite{mC}, 1.2.3, the proof is the same
as for the orthogonality relation for irreducible characters
of finite groups).
By normalizing them,
\[ d:=s {\bar s}^t, \quad S_{ij} := \frac{s_{ij}}{\sqrt{d_{ii}}},\]
where $\sqrt{d_{ii}}$ is the positive root,
we get an orthonormal matrix $S$ which we call
\defst{$S$-matrix} or \defst{Fourier matrix}
of $(R,B)$. We can recover the structure constants of $R$
(and the involution $\sim$) from
the $S$-matrix via the formula of Verlinde:
\begin{equation}\label{strcBV}
N_{ij}^l = \sum_k \frac{S_{ki} S_{kj} \BAR{S_{kl}}}{S_{k 0}},
\end{equation}
(this follows immediately by transporting the multiplication
via $\varphi$ from $\Cz^n$ to $R$).
The columns of the $s$-matrix are the image of $B$ in $\Cz^n$ under $\varphi$.
Conversely, given a matrix $S$ and a column (here column $0$),
we may define quantities $N_{ij}^l$ via (\ref{strcBV}),
which will be structure constants of a $\Zz$-based ring if
$S$ satisfies certain properties.

This is equivalent to the following construction:
start with a matrix $S\in\Cz^{n\times n}$
with $S\bar S^t=1$ and choose a column
$i_0$ in which all entries are non zero. Divide each row
by the entry of column $i_0$ to get a matrix $s$. The columns
of $s$ span a $\Zz$-lattice in $\Cz^n$ which is free since
$S$ is invertible. If this lattice
is closed under componentwise multiplication, then it
is a $\Zz$-algebra $R$ with the columns of $s$ as a basis.
In this case, we say that the matrix $S$ (or $s$) with unit $i_0$
define the $\Zz$-algebra $R$.
The involution $\sim$ corresponds to complex conjugation on the
columns of $S$.

\begin{exm}
Let $G$ be a finite group. An $s$-matrix of the character ring
of $G$ is the transposed character table of $G$.
\end{exm}

\begin{exm}
Untwisted affine Kac-Moody algebras have for each level $k$ a Kac-Peterson
matrix which is the $S$-matrix of a $\Zz$-based ring with
non negative structure constants (see \cite{vK}, 13.8).
\end{exm}

We will need the following lemma later on.
\begin{lem}\label{IsZbased}
Let $S\in\Cz^{n\times n}$ with $S\bar S^t=1$ such that
$N_{ij}^l = \sum_k \frac{S_{ki} S_{kj} \BAR{S_{kl}}}{S_{k 0}}\in\Zz$
and $S_{i0}\in\Rz^\times$ for all $i,j,l$.
If the set of columns of $S$ is invariant under complex conjugation,
then $S$ defines a $\Zz$-based ring.
\end{lem}
\begin{proof}
Let $\sim$ be the permutation of the columns given by complex
conjugation. Then
\[ N_{ij}^0=\sum_k \frac{S_{ki}S_{kj}\BAR{S_{k0}}}{S_{k0}}=
\sum_k S_{ki}\BAR{S_{k\tilde j}}=\delta_{i,\tilde j}\]
because $S_{k0}\in\Rz$ and $S\bar S^t=1$. We have to
check that $\sim$ is multiplicative:
\[ N_{\tilde i \tilde j}^m =
\sum_k \frac{S_{k\tilde i} S_{k\tilde j}\BAR{S_{k m}}}{S_{k0}} =
\sum_k \frac{\BAR{S_{k\tilde i}} \BAR{S_{k\tilde j}}S_{k m}}{\BAR{S_{k0}}} =
\sum_k \frac{S_{k i} S_{k j} \BAR{S_{k \tilde m}}}{S_{k0}} = N_{ij}^{\tilde m} \]
because $N_{\tilde i \tilde j}^m,S_{k0}\in\Rz$.
\end{proof}

\section{Fourier matrices for imprimitive complex reflection groups}

Let us define the Fourier matrices for the imprimitive
complex reflection groups $G(e,1,n)$
(compare with \cite{MUG}). The original definition is slightly technical, but
it takes a simple form if we express it by means of exterior powers
of the $S$-matrix of a group ring of a cyclic group.

\subsection{Definition of the Fourier matrices}

We use the notation of \cite{MUG}, 4A.
Let $e \ge 1$ and $Y$ be a totally ordered set with $d$ elements.
Consider the set
\[ \Psi:=\{ \psi : Y\rightarrow \{0,\ldots,e-1 \} \}\]
and a map $\pi : Y \rightarrow \Nz$. In \cite{MUG}, `$\Psi$' is a subset
of our $\Psi$; we will restrict to that subset later.
We define an equivalence relation $\sim_\pi$ on $\Psi$:
\[ \phi \sim_\pi \psi \quad \Longleftrightarrow \quad
\pi(\phi^{-1}(i))=\pi(\psi^{-1}(i)) \quad \mbox{for all}\quad 0\le i <e \]
for $\phi,\psi \in \Psi$, $\psi^{-1}(i):=\{ y\in Y\mid \psi(y)=i \}$,
and denote the class of $\psi$ by $[\psi]$.
Now call an element $\psi\in\Psi$ \defst{$\pi$-admissible} if for all $y,y'\in Y$ with
$\pi(y)=\pi(y')$ and $\psi(y)=\psi(y')$ we have $y=y'$.

A $\pi$-admissible $\psi$ can be interpreted in the sense
of \cite{MUG} as an $e$-symbol with entries in $\pi(Y)$:
an $e$-symbol is an ordered sequence ${\mathcal S}=(L_0,\ldots,L_{e-1})$
of $e$ strictly increasing finite sequences of natural numbers
$L_i=(\lambda_{i,1},\ldots,\lambda_{i,m_i})$, written
\[ {\mathcal S} = \begin{pmatrix}
\lambda_{0,1} & \cdots & \lambda_{0,m_0} \\
\lambda_{1,1} & \cdots & \lambda_{1,m_1} \\
\vdots & & \vdots \\
\lambda_{e-1,1} & \cdots & \lambda_{e-1,m_{e-1}}
\end{pmatrix} .\]
For $0\le i<e$, the set of entries of $L_i$ is $\pi(\psi^{-1}(i))$.

We define a matrix $\T$ indexed by the classes of $\pi$-admissible
elements of $\Psi$ (compare with \cite{MUG}, 4.10):
\[ \T_{[\phi],[\psi]} := \frac{(-1)^{m(e-1)}}{\tau(e)^m}
\sum_{\vartheta\in[\phi]} \epsilon(\vartheta)\epsilon(\psi)\prod_{y\in Y} \zeta^{-\vartheta(y)\psi(y)} \]
where $\zeta=\exp(2\pi\ii/e)$, $m:=\lfloor\frac{d}{e}\rfloor\in\Zz$ and
\[ \epsilon(\psi) := (-1)^{|\{ (y,y')\in Y\times Y \mid y<y',
 \psi(y)<\psi(y') \}|}, \quad
\tau(e):= \prod_{i=0}^{e-1}\prod_{j=i+1}^{e-1} (\zeta^i-\zeta^j) .\]
The Fourier matrices of \cite{MUG} are submatrices of $\T$.
Let $r:=|\pi(Y)|$ and $w_1,\ldots,w_r\in\Nz$ be such that $\pi(Y)=\{w_1,\ldots,w_r\}$.
Then $n_i:=|\pi^{-1}(w_i)|=|\psi(\pi^{-1}(w_i))|$ if $\psi$ is $\pi$-admissible.
Remark that if $\vartheta,\phi\in\Psi$ are equivalent (and $\pi$-admissible) then there is
a permutation $\sigma\in\Sym(Y)$ such that $\phi=\vartheta\circ\sigma$ and $\pi\circ\sigma=\pi$. Then
$\epsilon(\vartheta)\epsilon(\phi)=\varepsilon_\sigma$ holds,
where $\varepsilon_\sigma$ is the sign of the permutation $\sigma$.
Using this we get
\[ \T_{[\phi],[\psi]} = \frac{(-1)^{m(e-1)}}{\tau(e)^m} \epsilon(\phi)\epsilon(\psi)
\sum_{\stackrel{\sigma\in\Sym(Y)}{\pi\circ\sigma=\pi}}
\varepsilon_\sigma \prod_{y\in Y} \zeta^{-\phi(\sigma(y))\psi(y)},\]
and by defining $c:=(-1)^{m(e-1)}\ii^{-\binom{e-1}{2}m}\sqrt{e}^{d-em}$
\[ \T_{[\phi],[\psi]} = c \epsilon(\phi)\epsilon(\psi)
\sum_{\stackrel{\sigma\in\Sym(Y)}{\pi\circ\sigma=\pi}}
\varepsilon_\sigma \prod_{y\in Y} \frac{1}{\sqrt{e}} \zeta^{-\phi(\sigma(y))\psi(y)} \]
because $\tau(e)=\ii^{\binom{e-1}{2}}\sqrt{e}^e$.
Call a $\pi$-admissible $\psi\in\Psi$ \defst{ordered}, if for all $y,y'\in Y$, $y<y'$
with $\pi(y)=\pi(y')$ we have $\psi(y)<\psi(y')$. Then each class $[\psi]$
has exactly one ordered representative. So the
set indexing $\T$ is in bijection with
\[ \Xi:=\{ \psi \mid \psi\in\Psi, \:\: \psi \:\: \pi\mbox{-admissible and ordered} \} \]
and we will only consider elements from $\Xi$ from now on.
Take $\psi_1,\psi_2\in\Xi$ and let
$\psi_1(\pi^{-1}(w_\mu))=\{i_1^\mu,\ldots,i_{n_\mu}^\mu\}$,
$\psi_2(\pi^{-1}(w_\mu))=\{j_1^\mu,\ldots,j_{n_\mu}^\mu\}$ such that
$i_1^\mu<\ldots<i_{n_\mu}^\mu$, $j_1^\mu<\ldots<j_{n_\mu}^\mu$.
We finally get
\begin{equation}\label{finT} \BAR{\T_{\psi_1,\psi_2}} =
\bar c \epsilon(\psi_1)\epsilon(\psi_2)
\prod_{\mu=1}^r \sum_{\sigma\in S_{n_\mu}} \varepsilon_\sigma \prod_{\nu=1}^{n_\mu}
\frac{1}{\sqrt{e}} \zeta^{i_\nu^\mu j_{\sigma(\nu)}^\mu}.\end{equation}

\subsection{Connection to exterior powers}

Let $S\in\Cz^{e\times e}$ be the $S$-matrix of the group
ring of the cyclic group $\Zz/e\Zz$, so $S=(\frac{\zeta^{ij}}{\sqrt{e}})_{i,j}$.
Denote by $\Lambda^n \Cz^e$, $n\le e$, the subspace
of $\bigotimes_{i=1}^n \Cz^e$ spanned by
\[ C_n:=\Big\{\sum_{\sigma \in S_n} \varepsilon_\sigma
e_{i_{\sigma(1)}} \otimes \ldots \otimes e_{i_{\sigma(n)}} \mid 0\le i_1 < \ldots < i_n \le e-1\Big\},\]
where $e_0,\ldots,e_{e-1}$ is the
canonical basis of $\Cz^e$. The basis $C_n$ is indexed by the set of
$n$-tuples $(i_1,\ldots,i_n)$ with $0\le i_1<\ldots<i_n<e$; we will therefore
write $\bar i:=(i_1,\ldots,i_n)$ for the corresponding element of the basis.
The restriction of $\bigotimes_{i=1}^n S$ to $\Lambda^n \Cz^e$ defines an automorphism
corresponding to the matrix
\begin{equation}\label{Sdet}
 (\Lambda^n S)_{\bar i, \bar j} = \sum_{\sigma \in S_n} \varepsilon_\sigma
\prod_{\nu=1}^n S_{i_\nu,j_{\sigma(\nu)}}
= \det( (S_{i_\nu,j_{\nu'}})_{1\le \nu,\nu'\le n } )
\end{equation}
with respect to $C_n$.
Now consider the matrix $\Lambda^{n_1} S\otimes\cdots\otimes \Lambda^{n_r}S$
on the space $\Lambda^{n_1}\Cz^e\otimes\cdots\otimes\Lambda^{n_r}\Cz^e$ with the basis
\[ E:=\{ \bar i^1\otimes\ldots\otimes\bar i^r \mid \bar i^\mu \in C_{n_\mu},\:
1\le \mu\le r \}. \]
We identify $E$ with the set $\Xi$ by
\[ E\rightarrow \Xi, \quad \bar i^1\otimes\ldots\otimes\bar i^r \mapsto
\psi,\]
where $\psi$ is the element of $\Xi$ with
$\psi(\pi^{-1}(w_\mu))=\{i_1^\mu,\ldots,i_{n_\mu}^\mu\}$ for $1\le \mu\le r$
(remember that $\pi(Y)=\{w_1,\ldots,w_r\}$).
Formula (\ref{finT}) becomes
\[ \BAR{\T_{\psi_1,\psi_2}} = \bar c \epsilon(\psi_1)\epsilon(\psi_2)
(\Lambda^{n_1} S\otimes\cdots\otimes \Lambda^{n_r}S)_{\psi_1,\psi_2},\]
which explains why we will first concentrate on exterior powers
of $S$ to find out what the structure constants of submatrices of
$\T$ look like.
Note that since we are only interested in the integrality of structure
constants, we can multiply the column of $\psi_2$ by $\epsilon(\psi_2)$
for each $\psi_2$ and don't have to care about these signs anymore.
The $\epsilon(\psi_1)$ in the row of $\psi_1$ has no effect
on the structure constants because it is canceled in the formula
of Verlinde.

\section{Exterior powers}\label{extpow}

Consider again the matrix $\Lambda^n S$
where $S=(\frac{\zeta^{ij}}{\sqrt{e}})_{i,j}$ and $e,n\in\Nz$, $e\ge n$.
This matrix represents the restriction of $\bigotimes^n S$ to
$\Lambda^n\Cz^e$ with respect to the basis $C_n$ defined above.
Our goal is to prove that $\Lambda^n S$ defines a $\Zz$-algebra which
is a $\Zz$-based ring for an adequate basis (see Theorem \ref{extprodZRmB}),
so first we need to see that for suitable $\BAR{i_0}$ the structure constants
\[ N_{\bar j,\bar m}^{\bar k} = \sum_{\bar i \in C_n}
\frac{(\Lambda^n S)_{\bar i, \bar j} (\Lambda^n S)_{\bar i, \bar m} \BAR{(\Lambda^n S)_{\bar i, \bar k}}}
{(\Lambda^n S)_{\bar i, \bar i_0}}, \]
$\bar j,\bar m,\bar k \in C_n$, given by the formula of Verlinde are integers.
We take $\bar i_0:=(0,\ldots,n-1)$. Notice that
the above formula is then well defined because $(\Lambda^n S)_{\bar i, \bar i_0}$
in the denominator is a Vandermonde determinant and thus unequal to $0$.

\subsection{Connection to Schur functions}

We begin by analysing the
quotient
\[D_{\bar i,\bar j}:=\frac{(\Lambda^n S)_{\bar i, \bar j}}
{(\Lambda^n S)_{\bar i, \bar i_0}}.\]
The theorem about Jacobi-Trudi determinants (see
\cite{beS}, Theorem $4.5.1$) says
(\cite{beS}, Lemma $4.6.1$ and corollary $4.6.2$)
that $D_{\bar i,\bar j}$ is the Schur function
$s_{\bar j'}(\bar x) \in \Cz[[\bar x]]$,
$\bar x=\{x_1,x_2,\ldots\}$ evaluated at
$x_1=\zeta^{i_1},\ldots,x_n=\zeta^{i_n},x_{n+1}=0,\ldots$, where
$\bar j'$ is the partition
$\bar j':=(j_n-(n-1),\ldots,j_2-1,j_1)$.
There is an elementary proof of this
statement in \cite{mC}, 5.1.2.

The definition of the Schur function $s_{\bar j'}$
(see \cite{beS}, $4.4.1$) is
$s_{\bar j'} = \sum_{T\in \Tau_{\bar j'}} \bar x^T,$
where $\Tau_{\bar j'}$ is the set of semistandard
$\bar j'$-tableaux, that means tableaux of shape (Ferrer diagram)
$\bar j'$ with weakly increasing rows,
strictly increasing columns and entries in $\Nz$.
If $T_{h_1,h_2}$ are the entries of some $T\in \Tau_{\bar j'}$ then
$\bar x^T:= \prod_{h_1,h_2} x_{T_{h_1,h_2}}$. This
can also be written as
\[ \bar x^T = \prod_{\nu\in\Nz} x_\nu^{w_{\nu+1}(T)} \]
with suitable $w_{\nu+1}(T)\in \Nz$.
In our setting, $x_{n+1},x_{n+2},\ldots$ are all equal to
$0$, so from now on we consider
\[ s_{\bar j'} = \sum_{T\in\Tau_{\bar j'}}\prod_{\nu=1}^n x_\nu^{w_{\nu+1}(T)} \]
(by abuse of notation).
In the proof of the next theorem we will need the
following lemma.
\begin{lem}\label{Werte}
Let $a\in\Nz$ and
$ z_\nu(a):=\{ T\in\Tau_{\bar j'}\mid w_{\nu+1}(T)=a\} $
for $1\le \nu\le n$. Then
\[ |z_\nu(a)|=|z_{\nu'}(a)| \]
for all $1\le \nu' \le n$.
\end{lem}
\begin{proof}
Proposition $4.4.2$ in \cite{beS} says that $s_{\bar j'}(\bar x)$
is a symmetric function. This means, that
$s_{\bar j'} = s_{\bar j'}(x_{{\pi^{-1}}(1)},\ldots,x_{{\pi^{-1}}(n)})$
for all permutations $\pi\in S_n$ and hence
\[ \sum_{T\in\Tau_{\bar j'}}\prod_{\nu=1}^n x_\nu^{w_{\pi(\nu)+1}(T)} =
\sum_{T\in\Tau_{\bar j'}}\prod_{\nu=1}^n x_\nu^{w_{\nu+1}(T)}. \]
Therefore, for all $1\le \nu\le n$, $a\in\Nz$ and $T\in\Tau_{\bar j'}$
with $w_{\nu+1}(T)=a$ there is a $T'\in\Tau_{\bar j'}$ with $w_{\pi(\nu)+1}(T')=w_{\nu+1}(T)$.
But then we have $|z_\nu(a)|=|z_{\nu'}(a)|$ for all $1\le \nu'\le n$.
\end{proof}

\subsection{The structure constants are integers}

For $\bar j \in C_n$, we will write $\Tau_{\bar j}$ instead
of $\Tau_{\bar j'}$.
Here is the main theorem:

\begin{thm}\label{DSn}
The structure constants $N_{\bar j,\bar m}^{\bar k}$ defined above
are integers.
\end{thm}

\begin{proof}
First notice that by equation (\ref{Sdet}) we have
\[ (\Lambda^n S)_{\bar i, \bar j} = \frac{1}{\sqrt{e}^n} \sum_{\sigma \in S_n} \varepsilon_\sigma
\prod_{\nu=1}^n \zeta^{i_\nu j_{\sigma(\nu)}}, \]
so we get
\[ N_{\bar j,\bar m}^{\bar k} = \frac{1}{e^n}
\sum_{0\le i_1<\ldots<i_n\le e-1} \frac{P_{\bar i,\bar j} \BAR{P_{\bar i,\bar k}} P_{\bar i, \bar m}}
{P_{\bar i, \bar i_0}}, \]
if $P_{\bar i, \bar j} := \det(\zeta^{i_\mu j_\nu})_{\mu,\nu}$,
$\bar i_0:=(0,\ldots,n-1)$.
We have seen that $D_{\bar i,\bar m} := \frac{P_{\bar i,\bar m}}{P_{\bar i, \bar i_0}}$
can be expressed as
\[ D_{\bar i,\bar m} = \sum_{T \in \Tau_{\bar m}} \prod_{\nu=1}^n \zeta^{i_\nu w_{\nu+1}(T)}. \]
This remains well defined if we take any
tuple $(i_1,\ldots,i_n)$, $0\le i_1,\ldots,i_n \le e-1$ instead of restricting
to those with $i_1<\ldots<i_n$.
In general, $D_{\bar i, \bar m}$ can be non zero for
some $\bar i$ with two equal entries.
But the term $P_{\bar i,\bar j} \BAR{P_{\bar i,\bar k}} D_{\bar i, \bar m}$
is still $0$ in this case because the determinants $P_{\bar i,\bar j}, P_{\bar i,\bar k}$
vanish then.
Furthermore, we know that $D_{\bar i, \bar m}$ is invariant
under permutation of the $i_1,\ldots,i_n$, if $0\le i_1<\ldots<i_n\le e-1$.
Under permutation, the determinants
$P_{\bar i,\bar j}$ and $P_{\bar i,\bar k}$ are modified by
signs which cancel each other.
So we are allowed to write
\[ N_{\bar j,\bar m}^{\bar k} = \frac{1}{e^n} \sum_{0\le i_1<\ldots<i_n\le e-1}
P_{\bar i,\bar j} \BAR{P_{\bar i,\bar k}} D_{\bar i, \bar m} =
\frac{1}{e^n n!} \sum_{0\le i_1,\ldots,i_n\le e-1}
P_{\bar i,\bar j} \BAR{P_{\bar i,\bar k}} D_{\bar i, \bar m}. \]
We want to prove that $a:=\sum_{0\le i_1,\ldots,i_n\le e-1}
P_{\bar i,\bar j} \BAR{P_{\bar i,\bar k}} D_{\bar i, \bar m}$ is an integer
and congruent $0$ modulo $e^n n!$.
Substitute
\[ P_{\bar i,\bar j} = \sum_{\sigma \in S_n} \varepsilon_\sigma
\prod_{\nu=1}^n \zeta^{i_\nu j_{\sigma(\nu)}} \] to get
\[ a = \sum_{\sigma_1,\sigma_2 \in S_n} \varepsilon_{\sigma_1} \varepsilon_{\sigma_2}
\sum_{T \in \Tau_{\bar m}} \sum_{0\le i_1,\ldots,i_n\le e-1}
\prod_{\nu=1}^n \zeta^{i_\nu (j_{\sigma_1(\nu)}-
k_{\sigma_2(\nu)}+w_{\nu+1}(T))}, \]
where the inner sum can be rewritten as
\[ \sum_{i_1=0}^{e-1} \zeta^{i_1 (j_{\sigma_1(1)}-
k_{\sigma_2(1)}+w_{1+1}(T))} \sum_{i_2=0}^{e-1}
\zeta^{i_2 (j_{\sigma_1(2)}-k_{\sigma_2(2)}+w_{2+1}(T))} \cdots \hspace{1in}\]
\[\hspace{2in} \cdots \sum_{i_n=0}^{e-1}
\zeta^{i_n (j_{\sigma_1(n)}-k_{\sigma_2(n)}+w_{n+1}(T))} .\]
For a pair $(\sigma_1,\sigma_2)$, this is not zero if and only if
all brackets $(j_{\sigma_1(\nu)}-k_{\sigma_2(\nu)}+w_{\nu+1}(T))$
are congruent $0$ modulo $e$.
By Lemma \ref{Werte}, $w_{\nu+1}(T), T \in \Tau_{\bar m}$, take
the same values for all $\nu$ with the same multiplicities. Hence
if $(\sigma_1,\sigma_2)$ is an adequate pair (for which the sum is
not zero) then $(\sigma_1 \tau,\sigma_2 \tau)$, $\tau\in S_n$, is also adequate
($\varepsilon_{\sigma_1 \tau}\varepsilon_{\sigma_2 \tau} = \varepsilon_{\sigma_1}
\varepsilon_{\sigma_2}$).
Every adequate pair gives a contribution of $e^n$ at the end. So $a$
is congruent $0$ modulo $e^n n!$.
Hence $N_{\bar j,\bar m}^{\bar k} \in \Zz$.
\end{proof}

\subsection{Negative structure constants}

Here is an example in which the ring defined by $\Lambda^n S$ has negative structure constants:
\begin{exm}\label{extprodbsp1}
Take $e=4$, $n=2$. We have $6$ elements in the basis.
We write the multiplication table as the list of matrices
$(N_{0,i}^j)_{i,j} ,\ldots, (N_{5,i}^j)_{i,j}$.
If we write `$.$' for `$0$', then it is
{ \tiny
\[\begin{bmatrix}
         1& .& .& .& .& . \\
         .& 1& .& .& .& . \\
         .& .& 1& .& .& . \\
         .& .& .& 1& .& . \\
         .& .& .& .& 1& . \\
         .& .& .& .& .& 1
\end{bmatrix},
\begin{bmatrix}
         .& 1& .& .& .& . \\
         .& .& 1& 1& .& . \\
         .& .& .& .& 1& . \\
         .& .& .& .& 1& . \\
         -1& .& .& .& .& 1 \\
         .& -1& .& .& .& .
\end{bmatrix},
\begin{bmatrix}
         .& .& 1& .& .& . \\
         .& .& .& .& 1& . \\
         .& .& .& .& .& 1 \\
         -1& .& .& .& .& . \\
         .& -1& .& .& .& . \\
         .& .& .& -1& .& .
\end{bmatrix},
\] \[
\begin{bmatrix}
         .& .& .& 1& .& . \\
         .& .& .& .& 1& . \\
         -1& .& .& .& .& . \\
         .& .& .& .& .& 1 \\
         .& -1& .& .& .& . \\
         .& .& -1& .& .& .
\end{bmatrix},
\begin{bmatrix}
         .& .& .& .& 1& . \\
         -1& .& .& .& .& 1 \\
         .& -1& .& .& .& . \\
         .& -1& .& .& .& . \\
         .& .& -1& -1& .& . \\
         .& .& .& .& -1& . 
\end{bmatrix},
\begin{bmatrix}
         .& .& .& .& .& 1 \\
         .& -1& .& .& .& . \\
         .& .& .& -1& .& . \\
         .& .& -1& .& .& . \\
         .& .& .& .& -1& . \\
         1& .& .& .& .& .
\end{bmatrix}.\] }\noindent
This is not a $\Zz$-based ring, because there exists no involution $\sim$
as required. Applying substitutions $b \mapsto -b$ rectifies this.
But it is not possible to obtain a $\Zz$-based ring
with non negative structure constants just
by applying such substitutions (the computer easily checks
all $2^6$ sign changes).

It is unknown for which $e,n$ it is possible to get non
negative structure constants by applying sign changes.
Computations show that rings corresponding to $e,n$ with
at most 50 base elements have negative structure constants if
and only if both $e$ and $n$ are even and $1<n<e$.

As it is impossible to check all $2^{50}$ sign changes, we
apply another method:
If there is an appropriate sign change, then the new structure
constants will be the absolute values of the old ones and
define a $\Zz$-based ring.
So to decide if a given ring has such a sign change,
we try to compute an $s$-matrix for these new structure
constants (which fails if they do not define an algebra)
and then compare the $s$-matrices.
\end{exm}

A sign change as in the example above corresponds to
multiplying a column in the Fourier matrix by $-1$, which
does not change the $\Zz$-algebra.
As we see in the example, the matrix $\Lambda^n S$ does
not define a $\Zz$-based ring in general. However, we can
prove that there are sign changes such that we obtain
a $\Zz$-based ring.
We want to modify the matrix in such
a way that we can apply Lemma \ref{IsZbased}.
\begin{thm}\label{extprodZRmB}
Let $\Lambda^n S$ be as above where
$S=(\frac{\zeta^{ij}}{\sqrt{e}})_{i,j}$ and $e,n\in\Nz$, $e\ge n$.
Then $\Lambda^n S$ with the column $\bar i_0$ as unit
defines a $\Zz$-algebra $R$ and a basis $B$. Applying suitable
sign changes to $B$, we get a basis $B'$ such that
$(R,B')$ is a $\Zz$-based ring.
\end{thm}
\begin{proof}
First we define the involution $\sim$. Let $\bar i=(i_1,\ldots,i_n)\in C_n$
be an element of the basis. Define
$\bar i':=(n-1-i_1,\ldots,n-1-i_n)$ with entries taken modulo $e$.
Permuting $\bar i'$ we get an element $\tilde {\bar i}\in C_n$
and we will denote the sign of this permutation by $\gamma_{\bar i}$.
If as above $P_{\bar k, \bar i} := \det(\zeta^{k_\mu i_\nu})_{\mu,\nu}$ then
\[ \BAR{P_{\bar k,\bar i}}=P_{\bar k,-\bar i}
=P_{\bar k,\bar i'} \prod_{\nu=1}^n \zeta^{-k_\nu(n-1)}
=\gamma_{\bar i} P_{\bar k,\tilde{\bar i}} \prod_{\nu=1}^n \zeta^{-k_\nu(n-1)}. \]
With $\theta_{\bar k}:=\sqrt{\gamma_{\bar i_0}}\prod_{\nu=1}^n \zeta^{-\frac{k_\nu(n-1)}{2}}$
(for some choice of square root of $\gamma_{\bar i_0}$)
it follows that $\BAR{P_{\bar k,\bar i_0}\theta_{\bar k}}=P_{\bar k,\bar i_0}\theta_{\bar k}$ and
\[ \BAR{P_{\bar k,\bar i}\theta_{\bar k}}=\gamma_{\bar i} P_{\bar k,\tilde {\bar i}}
\BAR{\sqrt{\gamma_{\bar i_0}}}\prod_{\nu=1}^n \zeta^{-\frac{k_\nu(n-1)}{2}}=
\gamma_{\bar i} P_{\bar k,\tilde {\bar i}}\BAR{\sqrt{\gamma_{\bar i_0}}}
\sqrt{\gamma_{\bar i_0}}^{-1}\theta_{\bar k}=\gamma_{\bar i}\gamma_{\bar i_0}
P_{\bar k,\tilde{\bar i}}\theta_{\bar k}. \]
Remember that $P_{\bar k, \bar i}$ are the entries of $\Lambda^n S$
up to a factor that is a real number.
So multiplying each row $\bar k$ by $\theta_{\bar k}$, which is a root of
unity, we get a matrix $M$ whose columns we will denote by $v_{\bar i}$.
This matrix satisfies $M\bar M^t=1$, all entries in
$v_{\bar i_0}$ are real and for every $v_{\bar i}$ either
$\BAR{v_{\bar i}}$ or $-\BAR{v_{\bar i}}$ is a column of $M$.
The matrix $M$ defines the same algebra with the same basis
as $\Lambda^n S$ because we only have multiplied rows with roots of unity.

Now for each set $\{\bar i,\tilde{\bar i}\}$ such that
$\BAR{v_{\bar i}}=-v_{\tilde{\bar i}}$ choose
$\bar i$ or $\tilde{\bar i}$.
Multiply each column $v_{\bar i}$ of $M$ by $-1$
if $\bar i$ is a chosen element. The set of
columns in the resulting matrix is now closed under
complex conjugation. By Theorem \ref{DSn},
this matrix defines a ring with integer structure constants.
So all assumptions of Lemma \ref{IsZbased} are satisfied
and we obtain a $\Zz$-based ring.
\end{proof}

\section{The fusion algebras for the complex
reflection groups $G(e,1,n)$}\label{fusalgcrg}


The Fourier matrices for the complex reflection groups $G(e,1,n)$
decompose into blocks which are submatrices of the matrix $\T$
considered above. Here, $Y$ has $d=em+1$ elements and we restrict
to the subset
\[ E':=\{ \bar i^1\otimes\ldots\otimes \bar i^r\in E\mid
\sum_{y=1}^r\sum_{\nu=1}^{n_y} i_\nu^y \equiv m\binom{e}{2}
\:(\mbox{mod}\: e) \} \]
of $E$, so $\T':=(\T_{\xi_1,\xi_2})_{\xi_1,\xi_2\in E'}$ is
the matrix we will look at now.
Let $\abin:=m\binom{e}{2}$.
\begin{prop}\label{impganz}
Choose $a_y\in\Zz$, $1\le y\le r$, such that
\[ \xi_0 = (a_1,\ldots,a_1+n_1-1)\otimes\ldots\otimes(a_r,\ldots,a_r+n_r-1) \in E'.\]
Then the structure constants
\[ N_{\xi_1,\xi_2}^{\xi_3} := \sum_{\xi \in E'}
\frac{\T'_{\xi, \xi_1} \T'_{\xi, \xi_2}
\BAR{\T'_{\xi, \xi_3}}} {\T'_{\xi, \xi_0}}\]
are integers for all $\xi_1,\xi_2,\xi_3 \in E'$.
\end{prop}
\begin{proof}
Using the notation and the arguments of Theorem \ref{DSn}, we see that for
$\xi_1=\bar j^1\otimes\ldots\otimes\bar j^r$,
$\xi_2=\bar k^1\otimes\ldots\otimes\bar k^r$,
$\xi_3=\bar l^1\otimes\ldots\otimes\bar l^r \in E'$
\[ N_{\xi_1,\xi_2}^{\xi_3} =
e\sum_{\bar i^1\otimes\ldots\otimes\bar i^r \in E'}
\prod_{\mu=1}^r \frac{1}{e^{n_\mu}}
P_{\bar i^\mu,\bar j^\mu} \BAR{P_{\bar i^\mu,\bar k^\mu}} D_{\bar i^\mu, \bar l^\mu} \]
because $c\bar c= e^{d-em}=e$ and therefore
\[ N_{\xi_1,\xi_2}^{\xi_3} = \frac{e}{(n_r-1)!}
\left( \prod_{\mu=1}^{r-1} \frac{1}{n_\mu!} \right)
\sum_{0\le i^1_1,\ldots,i^1_{n_1}\le e-1} \cdots
\sum_{0\le i^{r-1}_1,\ldots,i^{r-1}_{n_{r-1}}\le e-1} \]
\[ \sum_{\substack{0\le i^r_2,\ldots,i^r_{n_r}\le e-1
\\ (i^r_1=\abin-i^r_2-\ldots-i^r_{n_r}
-\sum_{\mu=1}^{r-1}\sum_{\nu=1}^{n_\mu} i^\mu_\nu)}}
\prod_{\mu=1}^r \frac{1}{e^{n_\mu}}
P_{\bar i^\mu,\bar j^\mu} \BAR{P_{\bar i^\mu,\bar k^\mu}} D_{\bar i^\mu, \bar l^\mu}. \]
At the heart we find a power of $\zeta$ (see the proof of Theorem \ref{DSn})
with exponent of the form
${\sum_{\mu=1}^r \sum_{\nu=1}^{n_\mu} i_\nu^\mu \cdot w_{\mu,\nu}}$,
where the coefficient in front of $i_1^r$ equals
\[ w:=w_{r,1}=(j^r_{\sigma_1(1)}-k^r_{\sigma_2(1)}+w_{n+1}^r(T^r)).\]
Using the relation
\[ i^r_1 w= (\abin-i^r_2-\ldots-i^r_n-\sum_{\mu=1}^{r-1}\sum_{\nu=1}^{n_\mu} i^\mu_\nu) w \]
we can eliminate $i^r_1$ by subtracting $w$ from
each coefficient $w_{\mu,\nu}$ belonging to\\
$i^1_1,\ldots,i^{r-1}_{n_{r-1}},i^r_2,\ldots,i^r_{n_r}$;
the factor $\zeta^{a w} = {(\zeta^{\binom{e}{2}})}^{mw} = (\pm 1)^{mw}$
remains, which lies in $\Zz$ because
$2 \cdot \binom{e}{2}$ is divisible by $e$.

Then the proof goes on as in Theorem \ref{DSn}.
A pair $(\sigma_1,\sigma_2)$ may only be modified by
elements of the stabilizer of $1$ in $S_{n_r}$;
we obtain the desired factor $(n_r-1)!$.
The last sum yields only $e^{n_r-1}$
because there is no sum indexed by $i^r_1$.
Together with the $e$ in front, this
cancels against the factor $\frac{1}{e^{n_r}}$.
\end{proof}

As we did for the exterior powers, we want to
see that this $\Zz$-algebra is a $\Zz$-based ring
for a suitable basis.
In order to be able to use Lemma \ref{IsZbased},
we need to prove that the rows of $\T'$ are orthogonal
(this is implicit in \cite{MUG}, 4A):

\begin{prop}\label{orthT}
We have $\T' \BAR{\T'}^t = I$.
\end{prop}

\begin{proof}
We want to prove
$\sum_{\xi\in E'} \T'_{\xi_1,\xi}\BAR{\T'_{\xi_2,\xi}}
=\delta_{\xi_1,\xi_2}$,
where $\xi=\bar i^1\otimes\ldots\otimes \bar i^r$ runs
through $E'$, so
\[ (*)\quad \sum_{\mu=1}^r \sum_{\nu=1}^{n_\mu} i_\nu^\mu \equiv a \:(\mbox{mod } e).\]
As in Proposition \ref{impganz}
we are allowed to sum over all $0\le i_1^\mu,\ldots,i_{n_\mu}^\mu<e$
instead of $0\le i_1^\mu<\ldots< i_{n_\mu}^\mu<e$
for all $1\le \mu < r$.
This yields a factor $\frac{1}{n_\mu!}$.
The $i_\nu^\mu$ are only related by equation $(*)$.
We can therefore restrict to the case $r=1$
without loss of generality, which
simplifies the subsequent equations considerably.
Now, for $\xi_1:=\bar k:=(k_1,\ldots,k_n)$, $\xi_2:=\bar j:=(j_1,\ldots,j_n)$
and $\xi:=\bar i:=(i_1,\ldots,i_n)$ we have
\[ \sum_{\xi\in E'} \T'_{\xi_1,\xi}\BAR{\T'_{\xi_2,\xi}}
= \frac{e}{n!} \sum_{\substack{0\le i_1,\ldots,i_n\le e-1\\
\sum_{\nu=1}^n i_\nu \equiv a}}
\sum_{\sigma_1\in S_n} \sum_{\sigma_2\in S_n} \varepsilon_{\sigma_1 \sigma_2}
\prod_{\nu=1}^n \frac{1}{e} \zeta^{k_\nu i_{\sigma_1(\nu)}-j_\nu i_{\sigma_2(\nu)}}=\]
\[= \frac{e}{n!} \sum_{\sigma_1,\sigma_2\in S_n} \varepsilon_{\sigma_1\sigma_2}
\sum_{0\le i_2,\ldots,i_n\le e-1} \zeta^{a(k_{\sigma_1(1)}-j_{\sigma_2(1)})}
\frac{1}{e^n}\prod_{\nu=2}^n \zeta^{i_\nu(k_{\sigma_1(\nu)}-j_{\sigma_2(\nu)}
-k_{\sigma_1(1)}+j_{\sigma_2(1)})}, \]
which is not zero if and only if
\[ (**)\quad k_{\sigma_1(\nu)}-j_{\sigma_2(\nu)}\equiv k_{\sigma_1(1)}-j_{\sigma_2(1)}
\:(\mbox{mod } e) \]
for all $1\le \nu\le n$. But $0\le k_1<\ldots<k_n<e$ and
$0\le j_1<\ldots<j_n<e$, so this holds only if $\sigma:=\sigma_1=\sigma_2$.
If $(**)$ is satisfied, then because of $(*)$
\[ nf \equiv \sum_{\nu=1}^n k_\nu-j_\nu\equiv
\abin-\abin \equiv 0 \:(\mbox{mod } e),\]
where $f:=k_{\sigma(1)}-j_{\sigma(1)}$.
On the other hand, $n=em+1\equiv 1 \:(\mbox{mod } e)$ by assumption.
Consequently, the inner sum in the above formula
is zero only if
$\sigma_1=\sigma_2$ and $(**)$ are true, and in this case
$f=k_{\sigma_1(1)}-j_{\sigma_2(1)}\equiv 0$. Hence
\[ \sum_{\xi\in E'} \T'_{\xi_1,\xi}\BAR{\T'_{\xi_2,\xi}}
= \frac{e}{n!} \sum_{\sigma\in S_n} \zeta^{\abin f} \frac{1}{e}
= 1. \]
Conversely, from $(**)$, $d=0$ and $\sigma_1=\sigma_2$ follow
$\bar j = \bar k$.
\end{proof}
\begin{thm}\label{corMUG}
The Fourier matrices $\T'$ for the imprimitive complex reflection
group $G(e,1,n)$ define $\Zz$-algebras and bases which
are $\Zz$-based rings by applying sign changes to the bases.
\end{thm}
\begin{proof}
As in Proposition \ref{impganz}, we choose
\[ \xi_0 = (a_1,\ldots,a_1+n_1-1)\otimes\ldots\otimes(a_r,\ldots,a_r+n_r-1). \]
It remains to prove that a suitable involution $\sim$ exists.
We proceed exactly as in Theorem \ref{extprodZRmB}.
For $\xi=(i_1^1,\ldots,i_{n_1}^1)\otimes\ldots\otimes(i_1^r,\ldots,i_{n_r}^r) \in E'$,
let $\xi'$ be
\[ \xi':=(w_1-i_1^1,\ldots,w_1-i_{n_1}^1)\otimes\ldots\otimes
(w_r-i_1^r,\ldots,w_r-i_{n_r}^r),\]
with $w_\mu := n_\mu-1+2a_\mu$, $\mu=1,\ldots,r$.
Define $\tilde\xi$ to be the element of $E$ which
we get by sorting each bracket increasingly.
Then it lies in $E'$, because
\[ \sum_{\mu=1}^r \sum_{\nu=1}^{n_\mu} (n_\mu-1+2a_\mu)-i_\nu^\mu
= \Big(\sum_\mu n_\mu(n_\mu-1)+2n_\mu a_\mu\Big)-\abin = \]
\[ = 2\Big(\sum_\mu(\sum_{\nu=0}^{n_\mu-1} \nu+a_\mu)\Big)-\abin
\equiv 2 \abin-\abin \:(\mbox{mod }e), \]
where the last congruence comes from $\xi_0\in E'$.\\
It is easy to check that for
$\xi:=\bar k^1\otimes\ldots\otimes \bar k^r$,
$\xi_1:=\bar i^1\otimes\ldots\otimes \bar i^r$
\[ \BAR{\T_{\xi,\xi_1}} = \T_{\xi,\tilde\xi_1} \prod_{\mu=1}^r \gamma_{{\bar i}^\mu}
\prod_{\nu=1}^{n_\mu} \zeta^{-k_\nu^\mu(n_\mu-1+2a_\mu)},\]
where $\gamma_{{\bar i}^\mu}$ is the sign of the
permutation which sorts the tuple of $\xi'_1$
belonging to ${\bar i}^\mu$ (as in Theorem \ref{extprodZRmB}).
From now on, the proof continues exactly as in
Theorem \ref{extprodZRmB}. Remark that we need
Proposition \ref{orthT} at the end.
\end{proof}

\subsection{Eigenvalues and representation of $\SL_2(\Zz)$}
As Lusztig did it for his non abelian Fourier matrix, Malle
also defines a matrix  $\TT$ of eigenvalues of Frobenius associated
to a Fourier matrix $\T'$ 
for the complex reflection groups $G(e,1,n)$
(see \cite{MUG}, 4B).
The matrices $\T'$ and $\TT$ define an $\SL_2(\Zz)$-representation.

For a given $S$-matrix $S$, we will call a diagonal matrix
$T$ such that $S$ and $T$ define a representation of
$\SL_2(\Zz)$, i.e.
\[ S^4=1, \quad (ST)^3=1, \quad [S^2,T]=1 \]
a \defst{$T$-matrix} associated to $S$.
The pair $(S,T)$ is then called \defst{modular datum}
(this is not exactly the usual definition:
we do not require that the structure constants
defined by Verlinde's formula are non negative, and
$S$ is not necessarily symmetric).
Remark that a $T$-matrix is in general not uniquely determined by $S$.

There exists also a $T$-matrix for the exterior power
$\Lambda^n S$ as above. It is obtained by taking
the exterior power of a $T$-matrix corresponding to $S$
(compare \cite{MUG}, 4B):

\begin{prop}
Let $S:=(\frac{\zeta^{ij}}{\sqrt{e}})_{i,j}$ be the
$S$-matrix of the group ring of $\Zz/e\Zz$,
where $e\in\Nz$ and $\zeta:=\exp(2\pi\ii/e)$, $\zeta_{24}:=\exp(2\pi\ii/24)$.
Then the diagonal matrix $T$ with
\[ T_{i,i} = \zeta_{24}^{e-1} \zeta^{\frac{i^2+ei}{2}} \]
for $0\le i<e$ is a $T$-matrix for $S$.
\end{prop}

\begin{proof}
The equation $[S^2,T]=1$ is satisfied since $T$ is diagonal
and $S^2$ is the permutation corresponding to complex conjugation
on the columns of $S$.
We have to verify $((ST)^3)_{i,j} = \delta_{i,j}$.
Define $t_i := T_{i,i}$.
Because
\[ \zeta^\frac{i^2+ei}{2} = \zeta^\frac{(i+e)^2+e(i+e)}{2}, \]
$t_i$ only depends on the class of $i$ mod $e$. This
allows us to substitute $k$ by $k-l-i$ in the following
equation:
\[ ((ST)^3)_{i,j} = \frac{1}{e\sqrt{e}}
\sum_{k,l=0}^{e-1} \zeta^{ik+kl+lj} t_k t_l t_j
= \sum_{k,l} \zeta^{i(k-l-i)+(k-l-i)l+lj} t_{k-l-i} t_l t_j = \]
\[ = \frac{1}{e\sqrt{e}} \zeta_{24}^{3(e-1)}
\sum_k \zeta^{\frac{1}{2}(k^2+ek)+\frac{1}{2}(j^2-i^2+ej-ei)}
\sum_l \zeta^{lj-li}.\]
The inner sum is equal to $0$ if $i\ne j$ and equal to $e$ if $i=j$.
For $i=j$ it remains to prove
\[ \sum_{k=0}^{e-1} \zeta^{\frac{1}{2}(k^2+ek)} = \zeta_{24}^{-3(e-1)} \sqrt{e}, \]
which is a corollary of a theorem of Gau{\ss} (or equation 4.11 in \cite{MUG}).
\end{proof}

\section{Kac-Peterson matrices and exterior powers}

The construction of Fourier matrices from exterior powers
also shows up in a different context:
let ${\frak g}(A)$ be an affine Kac-Moody algebra belonging to an
$n\times n$ generalized Cartan matrix $A$ of rank $l$, $\frak h$ its
Cartan subalgebra and
$\langle \:,\: \rangle : {\frak h}\times{\frak h^*} \rightarrow \Cz$
the corresponding pairing (we use the notation of \cite{vK}). Define
\[ P:=\{\lambda\in {\frak h^*} \mid \langle \lambda,\alpha_i^\vee\rangle \in \Zz,\quad
i=0,\ldots,n-1 \}, \]
\[ P_+:=\{\lambda\in P \mid \langle \lambda,\alpha_i^\vee\rangle \ge 0,\quad
i=0,\ldots,n-1 \}. \]

Now let ${\frak g}(A)$ be of arbitrary untwisted type $X_l^{(1)}$
or $A_{2l}^{(2)}$ (then $n=l+1$).
The {\it fundamental weights} $\Lambda_i \in P$,
$i=0,\ldots,l$ are given by the equations
\[ \langle \Lambda_i,\alpha_j^\vee\rangle = \delta_{ij}, \quad \langle \Lambda_i,d\rangle=0\]
for $j=0,\ldots,l$, where $d\in {\frak h^*}$ is given
by $\langle \alpha_i,d\rangle=\delta_{i,0}$.
The $\{ \alpha_0^\vee,\ldots,\alpha_l^\vee,d \}$ form
a basis of $\frak h$ and
$\{ \alpha_0,\ldots,\alpha_l,\Lambda_0 \}$ form a basis of
$\frak h^*$.
The fundamental weights $\BAR\Lambda_i$ of the
finite dimensional Lie algebra $\frak g^\circ$ satisfy
\[ \Lambda_i = \BAR \Lambda_i + a_i^\vee\Lambda_0, \]
($\BAR \Lambda_0=0$ because $a_0^\vee=1$; for a definition of $a_i^\vee$,
see \cite{vK}, 6.1).

For each positive integer $k$, let
$P_+^k\subseteq P_+$ be the finite set
\[ P_+^k:=\Big\{ \sum_{j=0}^l \lambda_j\Lambda_j\mid
\lambda_j\in\Zz, \lambda_j\ge 0, \sum_{j=0}^l a_j^\vee\lambda_j=k \Big\}. \]
Kac and Peterson defined a natural $\Cz$-representation of the group $\SL_2(\Zz)$
on the subspace spanned by the affine characters of ${\frak g}$ which are indexed
by $P_+^k$.
The image of $\tiny \begin{pmatrix} 0&-1\\ 1&0 \end{pmatrix}$ under
this representation is determined in Theorem $13.8$ of \cite{vK}.
It is the so-called {\it Kac-Peterson matrix}.
For affine algebras of type $X_l^{(1)}$ or $A_{2l}^{(2)}$, this matrix is
\begin{equation}\label{kpmat}
S_{\Lambda,\Lambda'} = c\sum_{w\in W^\circ} \det(w)
\exp\left({-\frac{2\pi \ii (\BAR\Lambda+\bar\rho\mid
 w(\BAR\Lambda'+\bar \rho))}{k+h^\vee}}\right),
\end{equation}
where $\Lambda,\Lambda'$ runs through $P_+^k$, $(\cdot\mid\cdot)$ is
the normalized bilinear form of chapter $6$ of \cite{vK} and
$W^\circ$ is the Weyl group of $\frak g^\circ$.
The constant $c$ is unimportant for us,
since we want to use the matrix in the formula of Verlinde (\ref{strcBV}).

Each of these matrices defines a based ring. A classification
of the matrices belonging to type $X_l^{(1)}$ up to isomorphism
was given by Gannon in \cite{tG}.
Here we prove that the matrices of type $A_1^{(1)}$ are connected
to those of type $C_l^{(1)}$ via exterior powers:

\begin{prop}\label{AundC}
Let $k,l\in\Nz$, $k\ge 1, l\ge 2$ and $S$ be the
Kac-Peterson matrix of type $A_1^{(1)}$ and level $(k+l-1)$.
The Kac-Peterson matrix of type $C_l^{(1)}$ and level $k$
is the exterior power $\Lambda^l S$ of $S$.
\end{prop}

\begin{proof}
We use the notation of \cite{vK}.
Let $A$ be the Cartan matrix of type $C^{(1)}_l$.
Choose $\bar a,\bar a^\vee$ elements of the kernel of $A$
respectively $A^T$, say $\bar a:=(2,\ldots,2,1,1)$, $\bar a^\vee:=(1,\ldots,1)$
(arrange the matrix in such a way that $\alpha_0$ is at the end).
Furthermore, we have $\kappa:=k+h^\vee=k+\sum_{i=0}^l a_i^\vee=k+l+1$,
$\bar\rho:=(1,\ldots,1)$,
\[ P_+^k=\{\bar\lambda+\bar\rho \mid
\bar\lambda\in\{0,\ldots,k\}^l, \: \sum_{i=1}^l a_i^\vee\lambda_i\le k \}\]
and $W^\circ$ the Weyl group of type $C_l$.
Let $D$ be the diagonal matrix with
$\frac{a_1^\vee}{a_1},\ldots,\frac{a_l^\vee}{a_l}$ on its diagonal. Then
equation (\ref{kpmat}) becomes
\[ S_{\mu,\nu} = c \sum_{w\in W^\circ} \det(w)
\exp\left({-\frac{2\pi \ii}{\kappa}\mu D w^T {A^{-1}}^T \nu^T}\right) = \]
\[ = c \sum_{w\in W^\circ} \det(w)  \exp\left({-\frac{2\pi \ii}{\kappa}\mu M_w \nu^T}\right) \]
where $\mu,\nu \in P_+^k$ and $M_w:=D w^T {A^{-1}}^T$.
Consider the set
\[ \tilde P:=\{ \tilde\mu:=(\mu_1,\mu_1+\mu_2,\ldots,\sum_{i=1}^l \mu_i) \mid \mu\in P_+^k\}.\]
The base change $\mu\mapsto\tilde\mu$ transforms the
$M_w$, $w\in W^\circ$, into monomial matrices with entries $\pm \frac{1}{2}$
(up to a factor, this is the base change $\alpha_i\mapsto v_i$ from \cite{vK}, 6.7).
The formula for $S$ is now
\[ S_{\tilde\mu,\tilde\nu} = c \sum_{\sigma\in S_l}\sum_{f\in\{\pm 1\}^l}
(\varepsilon_\sigma \prod_{i=1}^l f_i)
\zeta^{-\sum_{i=1}^l f_i \tilde\mu_i\tilde\nu_{\sigma(i)}}
= c \sum_{\sigma\in S_l} \sum_{f\in\{\pm 1\}^l}
\varepsilon_\sigma \prod_{i=1}^l f_i
\zeta^{-f_i \tilde\mu_i\tilde\nu_{\sigma(i)}} = \]
\[ = \tilde c \sum_{\sigma\in S_l} \varepsilon_\sigma \prod_{i=1}^l
(\zeta^{\tilde\mu_i\tilde\nu_{\sigma(i)}}-\zeta^{-\tilde\mu_i\tilde\nu_{\sigma(i)}}) ,\]
where $\zeta:=\exp({\frac{2\pi \ii}{2 \kappa}})$.
In the last term, we recognize the determinants of the $l\times l$-submatrices
of the Kac-Peterson matrix of
type $A_1^{(1)}$ and level $k+l-1=\kappa-2$.
\end{proof}

\nocite{MUG}
\nocite{Lus1}
\nocite{mGgM}
\nocite{mC}
\nocite{vK}
\nocite{beS}
\nocite{CGR}
\nocite{tG}

\bibliographystyle{amsplain}
\bibliography{references}

\end{document}